\documentclass[runningheads]{llncs}
\usepackage{amssymb}
\usepackage{url}
\def\sA{{\sf A}} % because plain A is an abelian surface
\def\Atwo{{{\cal A}_2}}
\def\C{{\bf C}}
\def\Q{{\bf Q}}
\def\Qbar{{\,\overline{\!\Q\!}\,}}
\def\Fbar{{\,\overline{\!F}}}
\def\R{{\bf R}}
\def\Z{{\bf Z}}
\def\PP{{\bf P}}
\def\HH{{\cal H}}
\def\OO{{\cal O}}
\def\II{{\rm II}}
\def\MW{{Mordell--\kern-.12exWeil}}
\def\NeS{{N\'eron--Severi}}
\def\CI{{Clebsch--Igusa}}
\def\X{{\rm X}}
\def\XX{{\cal X}}
\def\Ness{{N_{\rm ess}}}
\def\Ndualess{{N^*_{\rm ess}}}
\def\ang#1{{\langle#1\rangle}}
\def\disc{\mathop{\rm{disc}}\nolimits}
\def\End{\mathop{\rm{End}}\nolimits}
\def\Gal{\mathop{\rm{Gal}}\nolimits}
\def\Km{\mathop{\rm{Km}}\nolimits}
\def\NS{\mathop{\rm{NS}}\nolimits}
\def\PSL{\mathop{\rm{PSL}}\nolimits}

\def\ra{{\rightarrow}}
\def\lra{{\leftrightarrow}}
\def\hra{{\hookrightarrow}}
\def\be{\begin{equation}}
\def\ee{\end{equation}}
\def\bea{\begin{eqnarray}}
\def\eea{\end{eqnarray}}
\def\0{^{\phantom0}}
\def\9{_{\phantom9}}
\def\sP{{s_{\kern-.05exP}\0}}
\begin{document}
\title{Shimura curve computations
   via K3 surfaces of \NeS\ rank at least $19$}
\titlerunning{Shimura curve computations via K3 surfaces}
\author{Noam D. Elkies}
\institute{Department of Mathematics, Harvard University,
Cambridge, MA 02138\\
\email{elkies@math.harvard.edu}\\
Supported in part by NSF grant DMS-0501029
}
\maketitle

\section{Introduction}

In \cite{NDE:Shim1} we introduced several computational challenges
concerning Shimura curves, and some techniques to partly address them.
The challenges are: obtain explicit equations for Shimura curves and
natural maps between them; determine a Schwarzian equation
on each curve (a.k.a.\ Picard--Fuchs equation, a linear second-order
differential equation with a basis of solutions whose ratio inverts
the quotient map from the upper half-plane to the curve);
and locate CM (complex multiplication) points on the curves.
We identified some curves, maps, and Schwarzian equations using
the maps' ramification behavior; located some CM points as images of
fixed points of involutions; and conjecturally computed others
by numerically solving the Schwarzian equations.

But these approaches are limited in several ways:
we must start with a Shimura curve with very few elliptic points
(not many more than the minimum of three); maps of high degree
are hard to recover from their ramification behavior,
limiting the range of provable CM coordinates; and these methods give
no access to the abelian varieties with quaternionic multiplication
(QM) parametrized by Shimura curves.
Other approaches somewhat extend the range where our challenges can be met.
Detailed theoretical knowledge of the arithmetic of Shimura curves
makes it possible to identify some such curves of genus at most~$2$
far beyond the range of~\cite{NDE:Shim1} (see e.g.~\cite{Roberts,GonRot}),
though not their Schwarzian equations or CM points.
Roberts~\cite{Roberts} showed in principle how to find CM coordinates
using product formulas analogous to those of~\cite{GZ} for
differences between CM \hbox{$j$-invariants},
but such formulas have yet to be used to verify and extend
the tables of~\cite{NDE:Shim1}.
Errthum~\cite{Errthum} recently used Borcherds products
to verify all the conjectural rational coordinates for CM points
tabulated in~\cite{NDE:Shim1} for the curves associated to
the quaternion algebras over~$\Q$ ramified at $\{2,3\}$ and $\{2,5\}$;
it is not yet clear how readily this technique might extend to
more complicated Shimura curves.  The \hbox{$p$-adic}
numerical techniques of~\cite{NDE:Shim2} give access to
further maps and CM points.  Finally, in the $\{2,3\}$ and $\{2,5\}$
cases Hashimoto and Murabayashi had already parametrized the relevant
QM abelian surfaces in 1995~\cite{HM}, but apparently such computations
have not been pushed further since then.

In this paper we introduce a new approach, which exploits the fact that some
Shimura curves also parametrize K3 surfaces of \NeS\ rank at~least~$19$.
``Singular'' K3 surfaces, those whose \NeS\ rank attains the
characteristic-zero maximum of~$20$, then correspond to
CM points on the curve.
We first encountered such parametrizations while searching for
elliptic K3 surfaces of maximal \MW\ rank over~$\Q(t)$
(see \cite{NDE:high_rank}), for which we used the K3 surface
corresponding to a rational \hbox{non-CM} point
on the Shimura curve $\X(6,79)/\ang{w_{6\cdot79}}$ of genus~$2$.
The feasibility of this computation suggested that such parametrizations
might be used systematically in Shimura curve computations.

This approach is limited to Shimura curves associated to
quaternion algebras over~$\Q$.  Within that important
special case, though, we can compute curves and CM points
that were previously far beyond reach.
The periods of the K3 surfaces should also allow the computation of
Schwarzian equations as in~\cite{LianYau},
though we have not attempted this yet.
We do, however, find the corresponding QM surfaces using
Kumar's recent formulas~\cite{Abhinav} that make explicit
Dolgachev's correspondence~\cite{Dolgachev}
between Jacobians of \hbox{genus-$2$} curves and
certain K3 surfaces of rank at least~$17$.
The parametrizations do get harder as the level of the Shimura curve
grows, but it is still much easier to parametrize the K3 surfaces
than to work directly with the QM abelian varieties --- apparently
because the level, reflected in the discriminant of the \NeS\ group,
is spread over $19$ \NeS\ generators rather than
the handful of generators of the endomorphism ring.\footnote{
  It would be interesting to quantify the computational complexity
  of such computations in terms of the level and the CM discriminant;
  we have not attempted such an analysis.
  }
In this paper we illustrate this with several examples of such
computations for the curves $\X(N,1)$ and their quotients.
As the example of $\X(6,79)/\ang{w_{6\cdot79}}$ shows,
the technique also applies to Shimura curves not covered by $\X(N,1)$,
but already for $\X(N,1)$ there is so much new data that we can
only offer a small sample here: the full set of results can be made
available online but is much too large for conventional publication.
Since we shall not work with $\X(N,M)$
for $M>1$, we abbreviate the usual notation $\X(N,1)$ to $\X(N)$ here.

The rest of this paper is organized as follows.
In the next section, we review the necessary background, drawn
mostly from \cite{Vigneras,Rotger:TAMS,BHPV}, concerning
Shimura curves, the abelian and K3 surfaces that they parametrize,
and the structure of elliptic K3 surfaces in characteristic zero;
then give A.~Kumar's explicit formulas for Dolgachev's correspondence,
% between Jacobians of \hbox{genus-$2$} curves and K3 surfaces
% whose \NeS\ lattice contains $U \oplus E_7 \oplus E_8$,
which we use to recover \CI\ coordinates for QM Jacobians from
our K3 parametrizations; and finally describe some of our techniques
for computing such parametrizations.  In the remaining sections
we illustrate these techniques in the four cases $N=6$, $N=14$,
$N=57$, and $N=206$.  For $N=6$ we find explicit elliptic models
for our family of K3 surfaces~$S$\/ parametrized by $\X(6)/\ang{w_6}$,
locate a few CM points to find the double cover
$\X(6) \ra \X(6) / \ang{w_6}$, transform $S$\/ to find
an elliptic model with essential lattice $\Ness \supset E_7 \oplus E_8$
to which we can apply Kumar's formulas,
and verify that our results are consistent with
previous computations of CM points~\cite{NDE:Shim1}
and \CI\ coordinates~\cite{HM}.
For $N=14$ we exhibit $S$ and verify the location of a CM point
that we computed numerically in~\cite{NDE:Shim1}
but could not prove using the techniques of~\cite{NDE:Shim1,NDE:Shim2}.
For $N=57$, the first case for which $\X(N)/\ang{w_N}$ has
positive genus, we exhibit the K3 surfaces parametrized by this curve,
and locate all its rational CM points.  For $N=206$, the last case
for which $\X(N)/\ang{w_N}$ has genus zero, we exhibit the corresponding
family of K3 surfaces and the hyperelliptic curves
$\X(206)$ and $\X(206)/\ang{w_2}$, $\X(206)/\ang{w_{103}}$
covering the rational curves
$\X(206)/\ang{w_{206}}$ and $\X(206)/\ang{w_2,w_{103}}$.

\section{Definitions and techniques}

{\bf Quaternion algebras over~$\Q$, Shimura curves, and QM abelian surfaces.}
Fix a squarefree integer $N>0$ with an even number of prime factors.
There is then a unique indefinite quaternion algebra $\sA/\Q$\/
whose finite ramified primes are precisely the factors of~$N$.
Let $\OO$ be a maximal order in~$\sA$.  Since $\sA$ is indefinite,
all maximal orders are conjugate in~$\sA$, and conjugate orders
will be equivalent for our purposes.
Let $\OO^*_1$ be the group
of units of reduced norm~$1$ in~$\OO$; let $\Gamma$ be
the arithmetic subgroup $\OO^*_1/\{\pm1\}$ of $\sA^*/\Q^*$;
and let $\Gamma^*$ be the normalizer of $\Gamma$
in the positive-norm subgroup of~$\sA^*/\Q^*$.
If $N=1$ then $\Gamma^* = \Gamma$; otherwise
$\Gamma^*/\Gamma$ is an abelian group of exponent~$2$,
and for each factor  $d|N$\/ there is a unique element
$w_d \in \Gamma^*/\Gamma$ whose lifts to $\sA^*$
have reduced norms in $d \cdot {\Q^*}^2$.

Because $\sA$ is indefinite, $\sA \otimes_\Q \R$
is isomorphic with the matrix algebra $M_2(\R)$,
so the positive-norm subgroup of $\sA^*/\Q^*$
is contained in $\PSL_2(\R)$ and acts on the upper half-plane~$\HH$.
The quotient $\HH/\Gamma$ is then a complex model of
the Shimura curve associated to~$\Gamma$, usually called $\X(N,1)$.
In~\cite{NDE:Shim1} we called this curve $\XX(1)$ in analogy with
the classical modular curve $\X(1)$ (see below),
since $N$\/ was fixed and we studied Shimura curves that we called
$\XX_0(p)$, $\XX_1(p)$, etc.,
associated with various congruence subgroups of $\sA^*/\Q^*$.
In this paper we restrict attention to $\HH/\Gamma$
and its quotients by subgroups of $\Gamma^*/\Gamma$;
thus we return to the usual notation, but simplify it to $\X(N)$
because we do not need $\X(N,M)$ for $M>1$.
If $N=1$ then $\sA \cong M_2(\Q)$, and we may take $\OO = M_2(\Z)$,
when $\Gamma = \Gamma^* = \PSL_2(\Z)$ and $\HH$\/ must be extended by
its rational cusps before we can identify $\HH/\Gamma$ with $\X(1)$.
Here we study curves $\X(N)$ and their quotients only for $N>1$,
and these curves have no cusps.

The Shimura curve $\X(N)$ associated to a quaternion algebra over~$\Q$
has a reasonably simple moduli description.
Fix a positive anti-involution $\varrho$ of~$\sA$ of the form
$\varrho(\beta) = \mu^{-1} \bar\beta \mu$ for some $\mu \in \OO$
with $\mu^2 + N = 0$.  Then $\X(N)$ parametrizes pairs $(A,\iota)$
where $A$\/ is a principally polarized abelian surface and
$\iota$ is an embedding of~$\OO$ into the ring $\End(A)$
of endomorphisms of~$A$, such that the Rosati involution
is given by~$\varrho$.  See \cite[\S2 and Prop.~4.1]{Rotger:TAMS}.
This gives $\X(N)$ the structure of an algebraic curve over~$\Q$.

An abelian surface with an action of a (not necessarily maximal)
order in a quaternion algebra is said to have
``quaternionic multiplication'' (QM).
A {\em complex multiplication} (CM) point of~$\X(N)$ is a point,
necessarily defined over~$\Qbar$,
for which $A$ has complex multiplication,
i.e.\ is isogenous with the square of a CM elliptic curve.
We shall use the QM abelian surfaces $A$ to find models for
the Shimura curves $\X(N)$ and locate some of their CM points.

When $N=1$,
an abelian surface together with an action of $\OO \cong M_2(\Z)$
is just the square of an elliptic curve,
so we recover the classical modular curve $\X(1)$.
We henceforth fix $N>1$.
Then the group $\Gamma^*/\Gamma$, acting on~$\X(N)$
by involutions that we also call $w_d$, is nontrivial.
These involutions are again defined over~$\Q$, taking
$(A,\iota)$ to $(A_d,\iota_d)$ for some $A_d$ isogenous with~$A$.
Specifically, $A_d$ is the quotient of~$A$ by the subgroup of the
\hbox{$d$-torsion} group $A[d]$ annihilated by the two-sided ideal
of~$\OO$ consisting of elements whose norm is divisible by~$d$,
and the principal polarization on~$A_d$ is $1/d$ times the
pull-back of the principal polarization on~$A$.
In particular $A_d$ is CM if and only if $A$ is.
Hence the notion of a CM point makes sense on the quotient of~$\X(N)$
by $\Gamma^*/\Gamma$ or by any subgroup of $\Gamma^*/\Gamma$.
If a CM point of discriminant~$-D$\/
on $\X(N) / (\Gamma^*/\Gamma)$ is rational
then the class group of $\Q(\kern-.1ex\sqrt{-D})$
must be generated by the classes of primes lying over factors $p|D$\/
that also divide~$N$.  Thus the class group has exponent~$1$ or~$2$
and bounded size; in particular, only finitely many~$D$\/ can arise.
In each of the cases $N=6$, $14$, $57$, and $206$ that we treat in this paper,
$N$\/ has two prime factors, so the class number is at most~$4$
and we can cite Arno \cite{Arno} to prove that a list of discriminants
of rational CM points is complete.
When $N$\/ has $4$ or $6$ prime factors we can use Watkins'
solution of the class number problem up to~$100$~\cite{Watkins}.

We have $A_N \cong A$ as principally polarized abelian surfaces,
but for $N>1$ the embeddings $\iota$, $\iota_N\0$ are not equivalent
for generic QM surfaces~$A$.
When we pass from $A$ to its Kummer surface
we shall lose the distinction between $\iota$ and $\iota_N\0$,
and so will at first obtain only the quotient curve $\X(N)/\ang{w_N}$.
We shall determine its double cover $\X(N)$ by locating
the branch points, which are the CM points on~$\X(N)/\ang{w_N}$
for which $A$ is isomorphic to the product of two elliptic curves
with CM by the quadratic imaginary order of discriminant~$-N$ or~$-4N$\/;
the arithmetic behavior of other CM points will then pin down
the cover, including the right quadratic twist over~$\Q$.

An abelian surface with QM by~$\OO$ has at least one
principal polarization, and the number of principal polarizations
of a generic surface with QM by~$\OO$
was computed in \cite[Theorem~1.4 and \S6]{Rotger:Ig}
in terms of the class number of $\Q(\kern-.1ex\sqrt{-N})$.
Each of these yields a map from $\X(N) / \ang{w_N}$ to
$\Atwo$, the moduli threefold of principally polarized abelian surfaces.
This map is either generically $1:1$ or generically $2:1$,
and in the $2:1$ case it factors through an involution
$w_d = w_{d'}$ on $\X(N)/\ang{w_N}$ where $d,d'>1$ are integers
such that $N = d d'$ and
\be
 \sA \cong \Bigl(\frac{-N,d}{\Q}\Bigr) \; [= \Bigl(\frac{d,d'}{\Q}\Bigr)
 ].
\label{eq:2:1}
\ee
(See the last paragraph of \cite[\S4]{Rotger:TAMS},
which also notes that a $2:1$ map occurs for $N=6$ and $N=10$,
each of which has a unique choice of polarization.  In the other
cases $N=14$, $57$, $206$ that we study in this paper,
only $1:1$ maps arise,
because the criterion (\ref{eq:2:1}) is not satisfied.)
We aim to determine at least one of the maps $\X(N)/\ang{w_N} \ra \Atwo$
in terms of the \CI\ coordinates on~$\Atwo$,
and thus to find the moduli of
the generic abelian surface with endomorphisms by~$\OO$.\footnote{
  Alas we cannot say simply ``find the generic abelian surface
  with endomorphisms by~$\OO$''\kern-.05ex, even up to quadratic twist,
  because there are abelian surfaces with rational moduli
  but no model over~$\Q$.
  }

{\bf K3 surfaces, elliptic K3 surfaces,
and the Dolgachev--Kumar correspondence.}
Let $F$\/ be a field of characteristic zero.
Recall that a {\em K3 surface} over~$F$\/
is a smooth, complete, simply connected algebraic surface $S/F$\/ 
with trivial canonical class.  The {\em \NeS\ group} $\NS(S)=\NS_\Fbar(S)$
is the group of divisors on~$S$\/ defined over the algebraic closure
$\Fbar$, modulo algebraic equivalence.  For a K3 surface this is
a free abelian group whose rank, the {\em Picard number}
$\rho = \rho(S)$, is in $\{1,2,3,\ldots,20\}$.
The intersection pairing gives $\NS(S)$ the structure of an
integral lattice; by the index theorem for surfaces, this lattice
has signature $(1, \rho - 1)$, and for a K3 surface the lattice is
{\em even}\/: $v \cdot v \equiv 0 \bmod 2$ for all $v \in \NS(S)$.
Over~$\C$, the cycle class map embeds $\NS(S)$ into the ``K3 lattice''
$H^2(S,\Z) \cong \II_{3,19} \cong U^3 \oplus E_8\ang{-1}^2$,
where $U = \II_{1,1}$ is the ``hyperbolic plane''
(the indefinite \hbox{rank-$2$} lattice with Gram matrix
$(\!{\scriptstyle{0\;1\atop1\;0}}\!)$),
and $E_8 \ang{-1}$ is the $E_8$ root lattice made negative-definite
by multiplying the inner product by~$-1$.
The Torelli theorem of Piateckii-Shapiro and \v{S}afarevi\v{c} \cite{PSS}
describes the moduli of K3 surfaces, at least over~$\C$: the embedding
of $\NS(S)$ into $\II_{3,19}$ is primitive, that is, realizes $\NS(S)$
as the intersection of $\II_{3,19}$ with a \hbox{$\Q$-vector} subspace
of $\II_{3,19} \otimes \Q$; for every such lattice~$L$\/ of signature
$(1,\rho-1)$, there is a nonempty (coarse) moduli space of pairs
$(S, \iota)$, where $\iota: L \rightarrow \NS(S)$ is a primitive
embedding consistent with the intersection pairing;
and each component of the moduli space has dimension $20-\rho$.
Moreover, for $\rho=20,19,18,17$ these moduli spaces repeat some
more familiar ones: isogenous pairs of CM elliptic curves for $\rho=20$,
elliptic and Shimura modular curves for $\rho=19$,
moduli of abelian surfaces with real multiplication or
isogenous to products of two elliptic curves for $\rho = 18$,
and moduli of abelian surfaces for certain cases of $\rho=17$.
Note the consequence that an algebraic family of
K3 surfaces in characteristic zero with $\rho \geq 19$
whose members are not all \hbox{$\Fbar$-isomorphic}
must have $\rho = 19$ generically, else there would be
a positive-dimensional family of K3 surfaces with $\rho \geq 20$.

An {\em elliptic} K3 surface $S/F$\/ is a K3 surface together with
a rational map $t: S \ra \PP^1$, defined over~$F$, whose generic fiber
is an elliptic curve.  The classes of the zero-section~$s_0$
and fiber~$f$\/ in $\NS(S)$ then satisfy $s_0 \cdot s_0 = -2$,
$s_0 \cdot f = 1$, and $f \cdot f = 0$,
and thus generate a copy of~$U$\/ in $\NS(S)$ defined over~$F$.
Conversely, {\em any} copy of~$U$\/ in $\NS(S)$ defined over~$F$\/
yields a model of~$S$\/ as an elliptic surface:
one of the standard isotropic generators or its negative is effective,
and has $2$ independent sections,
whose ratio gives the desired map to $\PP^1$.
We often use this construction to transform one elliptic model of~$S$\/
to another that would be harder to compute directly.
(Warning: in general one might have to subtract some base locus from
the effective generator to recover the fiber class~$f$.)

Since $\disc(U) = -1$ is invertible, we have
$\NS(S) = \ang{s_0,f} \oplus \ang{s_0,f}^\perp$,
with the orthogonal complement $\ang{s_0,f}^\perp$ having signature
$(0, \rho - 2)$; we thus write $\ang{s_0,f}^\perp = \Ness\ang{-1}$
for some positive-definite even lattice $\Ness$, the ``essential lattice''
of the elliptic K3 surface.  A vector $v \in \Ness$ of norm~$2$,
corresponding to $v \in \ang{s_0,f}^\perp$ with $v \cdot v = -2$,
is called a ``root'' of~$\Ness$; let $R \subseteq \Ness$
be the sublattice generated by the roots.  This root sublattice
decomposes uniquely as a direct sum of simple root lattices
$A_n$ ($n\geq 1$), $D_n$ ($n \geq 4$),
or $E_n$ \hbox{($6 \leq n \leq 8$)}.
These simple factors biject with reducible fibers,
each factor being the sublattice of~$\Ness$ generated by the components
of its reducible fiber that do not meet~$s_0$.
The graph whose vertices are these components,
and whose edges are their intersections,
is then the $A_n$, $D_n$, or $E_n$ root diagram;
if the identity component and its intersection(s)
are included in the graph then the extended root diagram
$\tilde{A}_n$, $\tilde{D}_n$, or $\tilde{E}_n$ results.
The quotient group $\Ness/R$\/
is isomorphic with the \MW\ group of the surface over~$\Fbar(t)$;
the isomorphism takes a point~$P$\/ to the projection of the
corresponding section $\sP$ to $\ang{s_0,f}^\perp$, and
the quadratic form on the \MW\ group induced from the pairing on~$\Ness$
is the canonical height.  Thus the \MW\ regulator is
$\tau^2 \disc(\Ness) \left/ \disc(R) \right.
= \tau^2 \left| \disc(\NS(S)) \right| / \disc(R)$,
where $\tau$\/ is the size of the torsion subgroup of the \MW\ group.

An elliptic surface has Weierstrass equation
$Y^2 = X^3 + A(t) X + B(t)$ for polynomials
$A,B$\/ of degrees at most~$8$, $12$
with no common factor of multiplicity at least~$4$ and~$6$ respectively,
and such that either $\deg(A)>4$ or $\deg(B)>6$ (i.e.,
such that the condition on common factors holds also at $t=\infty$
when $A,B$\/ are considered as bivariate homogeneous polynomials
of degrees $8$,~$12$).  The reducible fibers then occur at
multiple roots of the discriminant $\Delta = -16(4A^3+27B^2)$
where $B$\/ does not vanish to order exactly~$1$
(and at $t=\infty$ if $\deg\Delta \leq 22$ and $\deg B \neq 11$).
To obtain a smooth model for~$S$\/ we may start from the surface
$Y^2 = X^3 + A(t) X + B(t)$ in the $\PP^2$ bundle
\hbox{$\PP(O(0) \oplus O(2) \oplus O(3))$}
over~$\PP^1$ with coordinates $(1:X:Y)$,
and resolve the reducible fibers,
as exhibited in Tate's algorithm~\cite{Tate_alg},
which also gives the corresponding Kodaira types and simple root lattices.
This information can then be used to calculate
the canonical height on the \MW\ group, as in~\cite{Silverman:ht}.

The {\em Kummer surface} $\Km(A)$ of an abelian surface~$A$
is obtained by blowing up the $16 = 2^4$ double points of $A/\{\pm1\}$,
and is a K3 surface with Picard number $\rho(\Km(A)) = \rho(A) + 16 \geq 17$.
In general $\NS(\Km(A))$ need not consist of divisors defined over~$F$,
even when $\NS(A)$ does, because each \hbox{$2$-torsion} point of~$A$
yields a double point of $A/\{\pm1\}$ whose blow-up contributes to
$\NS(\Km(A))$, and typically $\Gal(\Fbar/F)$ acts nontrivially on $A[2]$.
But when $A$ is principally polarized Dolgachev~\cite{Dolgachev}
constructs another K3 surface $S_A/F$, related with $\Km(A)$ by
\hbox{degree-$2$} maps defined over~$\Fbar$,
together with a \hbox{rank-$17$} sublattice of $\NS(S_A)$
that is isomorphic with $U \oplus E_7 \oplus E_8$
and consists of divisor classes defined over~$F$.
It is these surfaces that we parametrize to get at the Shimura curves $\X(N)$.

If $A$ has QM then $\rho(A) \geq 3$, with equality for non-CM surfaces,
so $\rho(S_A) = \rho(\Km(A)) \geq 19$.
When $A$ has endomorphisms by~$\OO$, we obtain a sublattice
$L_N \subseteq \NS(S_A)$ of signature $(1,18)$ and discriminant~$2N$.
This even lattice $L_N$ is characterized by its signature
and discriminant together with the following condition:
for each odd $p|N$\/ the dual lattice $L_N^*$ contains a vector
of norm $c/p$ for some $c\in\Z$ such that $\chi_p(c) = -\chi_p(-2N/p)$,
where $\chi_p$ is the Legendre symbol $(\cdot/p)$; equivalently,
$\Ndualess$ contains a vector of norm $c/p$ with $\chi_p(c) = -\chi_p(+2N/p)$.
There is a corresponding local condition at~$2$, but
it holds automatically once the conditions at all odd $p|N$\/ are satisfied;
likewise when $N$\/ is odd it is enough to check all but one $p|N$.
The Shimura curve $\X(N)/\ang{w_N}$ parametrizes
pairs $(S,\iota)$ where $S$\/ is a K3 surface with $\rho(S) \geq 19$
and $\iota$ is an embedding $L_N \hra \NS(S)$.
If $\rho(S)=20$ then $(S,\iota)$ corresponds to a CM point
on~$\X(N)/\ang{w_N}$ whose discriminant equals $\disc(\NS(S))$.
The CM points of discriminant $-N$\/ or $-4N$\/ are the branch points
of the double cover $\X(N)$ of $\X(N)/\ang{w_N}$.  The arithmetic of
other CM points then determines the cover; for instance, if
$\X(N)/\ang{w_N}$ is rational, we know $\X(N)$ up to quadratic twist,
and then a rational CM point of discriminant $D \neq -N, -4N$\/
lifts to a pair conjugate over $\Q(\kern-.1ex\sqrt{-D})$.

The correspondence between $A$ and $S_A$ was made explicit by
Kumar~\cite[Theorem 5.2]{Abhinav}.
Let $A$ be the Jacobian of a \hbox{genus-$2$} curve~$C$,
and let $I_2, I_4, I_6, I_{10}$ be the \CI\ invariants of~$C$.
(If a principally polarized abelian surface~$A$ is not a Jacobian
then it is the product of two elliptic curves,
and thus cannot have QM unless it is a CM surface.)
We give an elliptic model of $S_A$ with $\Ness = R = E_7 \oplus E_8$,
using a coordinate~$t$ on~$\PP^1$ that puts the $E_7$ and $E_8$ fibers
at $t=0$ and $t=\infty$.  Any such surface has the formula
\be
Y^2 = X^3 + (a t^4 + a' t^3) X + (b'' t^7 + b t^6 + b' t^5)
\label{eq:E7E8}
\ee
for some $a,a',b,b',b''$ with $a', b'' \neq 0$.
(There are five parameters, but the moduli space has dimension
only $5-2 = 3$ as expected, because multiplying~$t$ by a nonzero scalar
yields an isomorphic surface, and multiplying $a,a'$ by~$\lambda^2$
and $b,b'$ by $\lambda^3$ for some $\lambda \neq 0$ yields
a quadratic twist with the same moduli.)
Kumar shows that setting
\be
(a,a',b,b',b'') =
\bigl(-I_4/12,\, -1,\, (I_2 I_4 - 3 I_6) / 108,\, I_2/24,\, I_{10}/4\bigr)
\label{eq:Kumar}
\ee
in~(\ref{eq:E7E8}) yields the surface $S_{J(C)}$.
Starting from any surface (\ref{eq:E7E8}) we may scale
$(t,X,Y)$ to $(-a't, {a'}^2 X, {a'}^3 Y)$ and divide through by ${a'}^6$
to obtain an equation of the same form with $a' = -1$; doing this
and solving (\ref{eq:Kumar}) for the \CI\ invariants $I_i$, we find
\be
(I_2, I_4, I_6, I_{10}) =
(-24 b'/a',\; -12 a,\; 96ab'/a' - 36b,\; -4a'b'').
\label{eq:rumaK}
\ee

If $A$ has QM by~$\OO$, but is not CM,
then the elliptic surface (\ref{eq:E7E8})
has a \MW\ group of rank~$2$ and regulator~$N$, with each choice
of polarization of~$A$\/ corresponding to a different \MW\ lattice.
The polarizations for which the map $\X(N)/\ang{w_N} \ra \Atwo$
factors through some $w_d$ are those for which the lattice
has an involution other than $-1$.
When this happens, two points on $\X(N)/\ang{w_N}$
related by $w_d$ yield the same surface (\ref{eq:E7E8})
but a different choice of \MW\ generators.
For example, when $N=6$ and $N=10$
these lattices have Gram matrices
$\frac12(\!${\large${\scriptstyle{5\;1\atop1\;5}}$}$\!)$ and
$\frac12(\!${\large${\scriptstyle{8\;0\atop0\;5}}$}$\!)$ respectively.
% LaTeX kludge!

{\bf Some computational tricks.}
Often we need elliptic surfaces with an $A_n$ fiber for
moderately large~$n$, that is, for which $4A^3+27B^2$ vanishes
to moderately large order $n+1$ at some $t=t_0$
at which neither $A$ nor~$B$\/ vanishes.
Thus we have approximately $(A,B) = (-3a^2,2a^3)$ near $t=t_0$.
Usually one lets $a$ be a polynomial
that locally approximates $(-A/3)^{1/2}$ at $t=t_0$, and writes
\be
(A,B) = (-3(a^2+2b),2(a^3+3ab)+c)
\label{eq:abc1}
\ee
for some $b,c$ of valuations $v(b)=\nu$, $v(c)=2\nu$ at $t_0$.
Then $v(\Delta) \geq 2\nu$ always,
and $v(\Delta) \geq 3\nu$ if and only if $v(3b^2-ac) \geq 3\nu$;
also if $\mu<\nu$ then
$v(\Delta) = 2\nu + \mu$ if and only if $v(3b^2-ac) = 2\nu+\mu$.
See \cite{Hall}; this was also the starting point of our analysis
in~\cite{NDE:Hall}.  For our purposes it is more convenient to allow
extended Weierstrass form and write the surface as
\be
Y^2 = X^3 + a(t) X^2 + 2b(t) X + c(t)
\label{eq:abc2}
\ee
with polynomials $a,b,c$ of degrees at most $4,8,12$ such that
$(v(b),v(c))=(\nu,2\nu)$.  Translating $X$\/ by $-a/3$ shows that
this is equivalent to~(\ref{eq:abc1}), with $a,b$\/ divided by~$3$
(so $\mu = v(b^2-ac)$ in~(\ref{eq:abc2})).  But (\ref{eq:abc2})
tends to produce simpler formulas, both for the surface itself
and for the components of the fiber, which are rational
if and only if $a$ is a square.  For instance, the Shioda--Hall
surface with an $A_{18}$ fiber \cite{Shioda_max,Hall}
can be written simply as
$$
Y^2 = X^3 + (t^4+3t^3+6t^2+7t+4)X^2 - 2(t^3+2t^2+3t+2)X + (t^2+t+1)
$$
with the $A_{18}$ fiber at infinity, and this is the quadratic twist
that makes all of $\NS(S)$ defined over~$\Q$.
The same applies to $D_n$, when $A':=A/t^2$ and $B':=B/t^3$
are polynomials such that $4{A'}^3+27{B'}^2$ has valuation $n-4$.
See for instance (\ref{eq:CM4}) below.
When we want singular fibers at several $t$ values we use an extended
Weierstrass form~(\ref{eq:abc2}) for which $(v(b),v(c))=(\nu,2\nu)$
holds (possibly with different~$\nu$) at each of these~$t$.

Having parametrized our elliptic surface~$S$\/
with $L_N \hra \NS(S)$,
we seek specializations of rank~$20$ to locate CM points.
In all but finitely many cases $S$\/ has an extra \MW\ generator.
In the exceptional cases, either some of the reducible fibers merge,
or one of those fibers becomes more singular,
or there is an extra $A_1$ fiber.  Such CM points are easy to locate,
though some mergers require renormalization to obtain a smooth model
and find the CM discriminant~$D$, as we shall see.
When there is an extra \MW\ generator, its height is at least $|D|/2N\!$,
but usually not much larger.  (Equality holds if and only if
the extra generator is orthogonal to the generic \MW\ lattice;
in particular this happens if $S$\/ has generic \MW\ rank zero.)
The larger the height of the extra generator,
the harder it typically is to find the surface.
This has the curious consequence that while the difficulty
of parametrizing $S$\/ increases with~$N$, the CM points
actually become easier to find.  In some cases we cannot
solve for the coefficients directly.
We thus adapt the methods of~\cite{NDE:Shim2},
exhaustively searching for a solution modulo a small prime~$p$
and then lifting it to a \hbox{$p$-adic} solution
to enough accuracy to recognize the underlying rational numbers.
We choose the smallest~$p$ such that $\chi_p(-D) = +1$,
so that reduction mod~$p$ does not raise the Picard number,
and we can save a factor of~$p$ in the exhaustive search
by first counting points mod~$p$ on each candidate~$S$\/
to identify the one with the correct~CM.

For large $N$\/ we use the following variation of the
\hbox{$p$-adic} lifting method to find the Shimura curve
$\X(N) / \ang{w_N}$ and the surfaces~$S$\/ parametrized by it.
First choose some indefinite primitive sublattice $L' \subset L_N$
and parametrize all $S$\/ with $\NS(S) \supseteq L'$.
Search in that family modulo a small prime~$p$ to find a surface $S_0$
with the desired~$L_N$.  Let $f_1,f_2$ be simple rational functions
on the $(S,L')$ moduli space.  We hope that the degrees, call them~$d_i$,
of the restriction of~$f_i$ to $\X(N) / \ang{w_N}$ are
positive but small; that $f_1$ is locally $1:1$ on the point of
$\X(N) / \ang{w_N}$ parametrizing $S_0$; and that the map
$(f_1,f_2) : \X(N) / \ang{w_N} \ra {\bf A}^2$
is generically $1:1$ to its image in the affine plane.
For various small lifts $\tilde{f}_1$ of $f_1(S_0)$ to~$\Q$,
lift $S_0$ to a surface $S/\Q_p$ with $f_1(S)=\tilde{f}_1$,
compute $f_2(S)$ to high \hbox{$p$-adic} precision,
and use lattice reduction to recognize $f_2(S)$ as the solution of
a polynomial equation $F(f_2) = 0$ of degree (at most) $d_1$.
Discard the few cases where the degree is not maximal,
and solve simultaneous linear equations to guess the coefficients
of~$F$\/ as polynomials of degree at most~$d_2$ in~$\tilde{f}_1$.
At this point we have a birational model $F(f_1,f_2)=0$
for $\X(N) / \ang{w_N}$.  Then recover a smooth model of the curve
(using Magma if necessary), recognize the remaining
coefficients of~$S$\/ as rational functions by solving
a few more linear equations, and verify that the surface
has the desired embedding $L_N \hra \NS(S)$.

\section{$N=6$: The first Shimura curve}

{\bf The K3 surfaces.}
We take $\Ness = R = A_2 \oplus D_7 \oplus E_8$,
which has discriminant $3 \cdot 4 \cdot 1 = 12 = 2N$,
and the correct behavior at~$3$ because
$A_2^*$ contains vectors of norm~$2/3$ with
$\chi_3\0(2) = -\chi_3\0(2\cdot 6/3) [= -1]$.
We choose the rational coordinate~$t$ on~$\PP^1$ such that
the $A_2$, $D_7$, and $E_8$ fibers are at $t=1$, $0$, and $\infty$
respectively.  If we relax the condition at $t=1$ by asking only that
the discriminant vanish to order at least~$2$ rather than~$3$
then the general such surface can be written as
\be
  Y^2 = X^3 + (a_0 + a_1 t) t X^2 + 2 a_0 b t^3 (t-1) X + a_0 b^2 t^5 (t-1)^2
\label{eq:S6_not}
\ee
for some $a_0, a_1, b$, with $a_1 b \neq 0$ lest the surface be too singular
at $t=0$.  The discriminant is then $t^9 (t-1)^2 \Delta_1(t)$
 with $\Delta_1$ a cubic polynomial such that
$\Delta(1) = -64 a_0 a_1 (a_0+a_1)^2 b^2$.
Thus $\Delta_1(1) = 0$ if and only if $a_1 = 0$ or $a_0 + a_1 = 0$.
In the latter case the surface has additive reduction at $t = 1$.
Hence we must have $a_1 = 0$.  The non-identity components of
the resulting $A_2$ fiber at $t=1$ then have $X = O(t-1)$;
we calculate that $X = x_1 (t-1) + O((t-1)^2)$ makes
$Y^2 = (x_1+b)^2 a_0 (t-1)^2 + O(t-1)^3$.
Therefore these components are rational if and only if $a_0$ is a square.
We can then replace $(X,Y,b)$ by $(a_0\0 X, a_0^{3/2} Y, a_0\0 b)$
in~(\ref{eq:S6_not}) to obtain the formula
\be
  Y^2  =  X^3  +  t X^2  +  2 b t^3 (t-1) X  +  b^2 t^5 (t-1)^2
\label{eq:S6}
\ee
for the general elliptic K3 surface with $\Ness = R = A_2 D_7 E_8$
and rational $A_2$ components.  The two components of the $D_7$ fiber
farthest from the identity component then have
$X = b t^2 + O(t^3)$, so $Y^2 = b^3 t^6 + O(t^7)$; thus these components
are both rational as well if and only $b$ is a square, say $b=r^2$.
Then $b$ and $r$ are rational coordinates on the Shimura curves
$\X(6)/\ang{w_2,w_3}$ and $\X(6)/\ang{w_6}$ respectively,
with the involution $w_2 = w_3$ on $\X(6)/\ang{w_6}$
taking $r$ to $-r$.

The elliptic surface (\ref{eq:S6}) has discriminant
$\Delta = 16 b^3 t^9 (t-1)^3 (27 b (t^2-t) - 4)$.
Thus the formula (\ref{eq:S6}) fails at $b=0$,
and also of course at $b=\infty$.
Near each of these two points we change variables to obtain
a formula that extends smoothly to $b=0$ or $b=\infty$ as well.  
These formulas require extracting respectively a fourth and third root
of~$\beta$, presumably because $b=0$ and $b=\infty$ are elliptic points
of the Shimura curve.
For small~$b$, we take $b = \beta^4$ and replace $(t,X,Y)$ by
$(t/\beta^2, X/\beta^2, Y/\beta^3)$ to obtain
\be
  Y^2  =  X^3  +  t X^2  +  2 t^3 (t-\beta^2) X  +  t^5 (t-\beta^2)^2,
\label{eq:S6_0}
\ee
with the $A_2$ fiber at $t=\beta^2$ rather than $t=1$.
When $\beta = 0$, this fiber merges with the $D_7$ fiber at $t=0$
to form a $D_{10}$ fiber, but we still have a K3 surface, namely
$Y^2 = X^3 + t X^2 + 2 t^4 X + t^7$, with $L = R = D_{10} \oplus E_8$.
This is the CM point of discriminant~$-4$.
For large~$b$, we write $b = 1/\beta^3$ and replace $(X,Y)$ by
$(X/\beta^2, Y/\beta^3)$ to obtain
\be
  Y^2  =  X^3  +  \beta^2 t X^2  +  2 \beta t^3 (t-1) X  +  t^5 (t-1)^2;
\label{eq:S6_inf}
\ee
then taking $\beta \ra 0$ yields the surface $Y^2 = X^3 + t^5 (t-1)^2$
with $\Ness = R_0 = A_2 \oplus E_8 \oplus E_8$:
the $t=0$ fiber changes from $D_7$ to $E_8$,
and the $t=1$ fiber becomes additive but still contributes $A_2$ to~$R$\/
(Kodaira type~IV rather than I$_3$).
This is the CM point of discriminant~$-3$.

{\bf Two more CM points.}
The factor $27 b (t^2-t) - 4$ of~$\Delta$
is a quadratic polynomial in~$t$ of discriminant $27b(27b+16)$.
Hence at $b=-16/27$ we have
$\Ness = R = A_1 \oplus A_2 \oplus D_7 \oplus E_8$,
and we have located the CM point of discriminant $-24$.
Three points fix a rational coordinate on~$\PP^1$,
so we can compare with the coordinate used in \cite[Table~1]{NDE:Shim1},
which puts the CM points of discriminant $-3$, $-4$, and $-24$
at $\infty$, $1$, and~$0$ respectively; thus that coordinate is
$1 + 27b/16$.  This also confirms that $\X(6)$ is obtained by
extracting a square root of $-(27r^2+16)$.

We next locate a CM point of discriminant $-19$ by finding $b$ for which
the surface (\ref{eq:S6}) has a section $\sP$ of canonical height $19/12$.
This is the smallest possible canonical height for a surface with
$R = A_2 \oplus D_7 \oplus E_8$,
because the na\"{\i}ve height is at least~$4$
and the height corrections at the $A_2$ and $D_7$ fibers
can reduce it by at most $2/3$ and $7/4$ respectively,
reaching $4 - 2/3 - 7/4 = 19/12$.
Let $(X(t),Y(t))$ be the coordinates of a point~$P$\/ of height $19/12$.
Then $X(t)$ and $Y(t)$ are polynomials of degree at most $4$ and~$6$
respectively (else $\sP$ intersects $s_0$
and the na\"{\i}ve height exceeds~$4$), and $X$\/ vanishes at $t=1$
(so $\sP$ passes through a non-identity component of the $A_2$ fiber)
and has the form $b t^2 + O(t^3)$ at $t=0$ (so $\sP$ meets one of
the components of the $D_7$ fiber farthest from the identity component).
That is, $X = b(t^2-t^3)(1 + t_1 t)$ for some $t_1$.
Substituting this into (\ref{eq:S6}) and dividing by
the known square factor $(t^4-t^3)^2$ yields $b^3$ times
\be
-t_1^3 t^4  +  (t_1^3 - 3 t_1^2) t^3
+  3 (t_1^2-t_1) t^2  +  ((3 t_1 - 1) + b^{-1} t_1^2) t  +  1,
\label{eq:19}
\ee
so we seek $b,t_1$ such that the quartic (\ref{eq:19}) is a square.
We expand its square root in a Taylor expansion about $t=0$
and set the $t^3$ and $t^4$ coefficients equal to zero.
This gives a pair of polynomial equations in $b$ and~$t_1$,
which we solve by taking a resultant with respect to $t_1$.
Eliminating a spurious multiple solution at $b=0$,
we finally obtain $(b,t_1) = (81/64,-9)$, and confirm that
this makes (\ref{eq:19}) a square, namely $(27t^2-18t-1)^2$.
Therefore $81/64$ is the \hbox{$b$-coordinate} of a CM point
of discriminant~$-19$.  Then $1 + 27b/16 = 3211/2^{10}$,
same as the value obtained in~\cite{NDE:Shim1}.

{\bf \CI\ coordinates.}
The next diagram shows the graph whose vertices are
the zero-section (circled) and components of reducible fibers
of an elliptic K3 surface~$S$\/ with $\Ness = A_2 \oplus D_7 \oplus E_8$,
and whose edges are intersections between pairs of these rational curves
on the surface.  Eight of the vertices form an extended root diagram
of type $\tilde{E}_7$, and are marked with their multiplicities
in a reducible fiber of type $E_7$ of an alternative elliptic model
for~$S$.  We may take either of the unmarked vertices of the
$\tilde{D}_7$ subgraph as the zero-section.  Then the essential lattice
of the new model includes an $E_8$ root diagram as well as the forced $E_7$.
We can thus apply Kumar's formulas to this model once we compute
its coefficients.

\centerline{
\setlength{\unitlength}{.45in}
\begin{picture}(10,3.7)(0,-0.3)
\thicklines
\put(1,2){\line(1,0){7}}
  \put(2.1,2.1){\makebox(0,0)[bl]{$1$}}
  \put(3.0,2.1){\makebox(0,0)[bl]{$2$}}
  \put(4.0,2.1){\makebox(0,0)[bl]{$3$}}
  \put(5.1,2.1){\makebox(0,0)[bl]{$4$}}
  \put(6.0,2.1){\makebox(0,0)[bl]{$3$}}
  \put(7.1,2.1){\makebox(0,0)[bl]{$2$}}
  \put(7.95,2.1){\makebox(0,0)[b]{$1$}}
  \put(5.1,3){\makebox(0,0)[l]{$2$}}
\multiput(1,2)(1,0){8}{\circle*{.1}}
\put(7,2){\circle{.2}}
\put(8,2){\line(4,-1){1}}
\put(8,2){\line(1,1){.75}}
\put(9,1.75){\line(-1,4){.25}}
  \put(8.75,2.75){\circle*{.1}}
  \put(9,1.75){\circle*{.1}}
\put(7,2){\line(0,-1){1}}
\put(0,1){\line(1,0){7}}
\multiput(0,1)(1,0){8}{\circle*{.1}}
\put(2,2){\line(0,1){1}}
  \put(2,3){\circle*{.1}}
\put(5,2){\line(0,1){1}}
  \put(5,3){\circle*{.1}}
\put(2,1){\line(0,-1){1}}
  \put(2,0){\circle*{.1}}
\thinlines
\put(-0.25,-0.25){\framebox(7.5,1.5){}}
  \put(0,0.05){\makebox(0,0)[l]{$\tilde{E}_8$}}
\put(0.75,1.75){\framebox(5.65,1.5){}}
  \put(1.05,2.95){\makebox(0,0)[c]{$\tilde{D}_7$}}
\put(7.75,1.5){\framebox(1.5,1.5){}}
  \put(7.85,2.9){\makebox(0,0)[tl]{$\tilde{A}_2$}}
\end{picture}
}
%          0           2             0
%          |           |            / \
%      0 - 1 - 2 - 3 - 4 - 3 - 2 - 1 - 0
%                              |
%  0 - 0 - 0 - 0 - 0 - 0 - 0 - 0
%          |
%          0
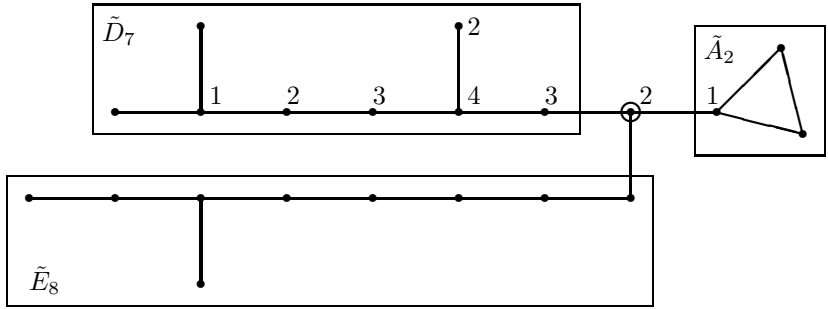
\centerline{Figure 1: An $\tilde{E}_7$ divisor supported on the zero-section and}
\centerline{fiber components of an $A_2 D_7 E_8$ surface}

\vspace*{2ex}

The sections of the $\tilde{E}_7$ divisor are generated by $1$ and
$u := X/(t^4-t^3) + b/t$.  Thus $u: S \ra \PP^1$ gives
the new elliptic fibration.  Taking $X = (t^3-t^2)(tu-b)$
in (\ref{eq:S6}) and dividing by $(t^4-t^3)^2$ yields
$Y_1^2 = Q(t)$ for some quartic~$Q$.  Using standard formulas
for the Jacobian of such a curve, and bringing the resulting
surface into Weierstrass form, we obtain a formula (\ref{eq:E7E8})
with $(a,a',b,b',b'')$ replaced by $(-3b, 1, -2b^2, -(b+1), -b^3)$.
As expected this surface has \MW\ rank~2
with generators of height $5/2$, namely
$$
  \bigl(
  r^6 t^4 + 2(r^4 + r^3) t^3 + (r^2+1)t^2,\,
  r^9 t^6 + 3(r^7 + r^6) t^5 + 3(r^5 + r^4 + r^3) t^4 + (r^3 + 1) t^3
  \bigr)
$$
and the image of this section under $r \, \lra {} -r$
(recall that $b=r^2$).
The formula~(\ref{eq:rumaK}) yields the \CI\ coordinates
\be
(I_2, I_4, I_6, I_{10}) = ((24b+1), 36b, 72b(5b+4), 4b^3).
\label{eq:CI6}
\ee

\section{$N=14$: The CM point of discriminant $-67$}

{\bf The K3 surfaces.}
Here we take $\Ness = R = A_3 \oplus A_6 \oplus E_8$,
which has discriminant $4 \cdot 7 \cdot 1 = 28 = 2N$,
and the correct behavior at~$7$ because $A_6^*$ contains
vectors of norm~$6/7$ with $\chi_7\0(6) = -\chi_7\0(2\cdot 14/7)  [= -1]$.
We put the $A_3$, $A_6$, and $E_8$ fibers at $t=1$, $0$, and $\infty$
respectively.  We then seek an extended Weierstrass form (\ref{eq:abc2})
with $a,b,c$ of degrees $2,4,7$ such that $t^3-t^2|b$, $(t^3-t^2)^2|c$,
and $(t^3-t^2)^6|b^2-ac$.  This gives at least $A_3, A_5, E_8$.
It is then easy to impose the extra condition $t^7 | \Delta$,
and we obtain
$a = \lambda ((s+1)t^2 + (3s^2+2s)t + s^3)$,
$b = \lambda^2 (s+1)  ((4s+2)t + 2s^2)  (t^3-t^2)$,
$c = \lambda^3 (s+1)^2 (t+s) (t^3-t^2)^2$
for some $s,\lambda$.  The twist $\lambda$ must be chosen so that
$a(0)$ and $a(1)$ are both squares; this is possible if and only if
$s^2+s$ is a square, so $s = r^2/(2r+1)$ for some~$r$.
Thus $r$ and $s$ are rational coordinates on
$\X(14) / \ang{w_{14}}$ and $\X(14) / \ang{w_2,w_7}$ respectively,
with the involution $w_2 = w_7$ on $\X(14) / \ang{w_{14}}$
taking $r$ to $-r/(2r+1)$.
The formula in terms of~$r$ is cleaner if we let the $A_3$ fiber
move from $t=1$; putting it at $t=2r+1$ yields
\bea
a &=& ((r+1)^2)t^2 + (3r^4+4r^3+2r^2)t + r^6),
\nonumber\\
b &=& 2(r+1)^2 ((2r^2+2r+1)t + r^4) (t-(2r+1)) t^2,
\label{eq:abc_14} \\
c &=& (r+1)^4 (t-(2r+1))^2 (t+r^2) t^4.
\nonumber
\eea

{\bf Easy CM points.}  At $r=0$, the $A_6$ fiber becomes $E_7$,
so we have a CM point with $D=-8$;
at $r=-1/2$, the $A_3$ and $A_6$ fibers merge to $A_{10}$,
giving a CM point with $D=-11$.  These have $s=0$, $s=\infty$
respectively.  There is an extra $A_1$ fiber when $11s^2+3s+8 = 0$;
the roots of this irreducible quadratic give the CM points with $D=-56$
(and their lifts to $\X(14) / \ang{w_{14}}$ are the branch points
of the double cover $\X(14)$).  In~\cite{NDE:Shim1} we gave
a rational coordinate~$t$ on $\X(14) / \ang{w_2,w_7}$ for which
the CM points of discriminants $-8$, $-11$, and $-56$ had
$t=0$, $t=-1$, and $16t^2+13t+8=0$ respectively.  Therefore
that~$t$ is our $-s/(s+1)$.

{\bf A harder CM point.}  At the CM point of discriminant $-67$
our surface has a section of height $67/28 = 4 - (3/4) - (6/7)$.
Thus $Y^2 = X^3 + aX^2 + bX + c$ has a solution in polynomials $X,Y$\/
of degrees $4,6$ with $X(0)=X(2h+1)=0$ and $Y$\/ having valuation
exactly~$1$ at $t=0$ and $t=2h+1$.  An exhaustive search mod~$17$
quickly finds an example, whose lift to $\Q_{17}$ then yields
$r = -35/44$ with
$$
 X = \frac{3^4}{5^2 \, 22^5} \, t \, (22t+13) \, (527076t^2 + 760364t + 275625).
$$
Thus $s=-1225/1144$, and $-s/(s+1)$ confirms the entry $-1225/81$
in the $|D|=67$ row of \cite[Table~5]{NDE:Shim1}.

\section{$N=57$: The first curve $\X(N)/w_N$ of positive genus}

{\bf The K3 surfaces.}
We cannot have $\Ness = R$\/ here because there is
no root lattice of rank~$17$ and discriminant $6 \cdot 19$.
Instead we take for $R$\/ the \hbox{rank-$16$} lattice
$A_5 \oplus A_{11}$ of discriminant $6 \cdot 12 = 72$,
and require an infinite cyclic \MW\ group $\Ness/R$\/
with a generator corresponding to a section that
meets the $A_5$ and $A_{11}$ fibers in non-identity components
farthest from the $A_5$ identity and nearest the $A_{11}$ identity
respectively, and does not meet the zero-section
(i.e., for which $X$\/ is a polynomial of degree at most~$4$ in~$t$).
Such a point has canonical height
\be
4 - \frac{3 \cdot 3}{6} - \frac{1 \cdot 11}{12}
= \frac{19}{12} = \frac{2N}{\disc R}.
\label{eq:19/12}
\ee
Thus $\disc(\Ness)$ has the desired discriminant~$2N$.
We may check the local conditions by noting that $A^*_5$
contains a vector of norm~$4/3$ that remains in~$\Ndualess$,
and $\chi_3\0(4) = -\chi_3\0(2\cdot 57/3)  [= +1]$.
We put the $A_5$ fiber at $t=0$ and the $A_{11}$ fiber at $t=\infty$.
We eventually obtain the following parametrization in terms of
a coordinate~$r$ on the rational curve $\X(57) / \ang{w_3,w_{19}}$:
let
\bea
p(r) &=& 4(r-1)(r^2-2) + 1,
\nonumber\\
d &=& (r^2-1)^2 (9t+(2r-1)p(r)),
\nonumber\\
c &=& 9t^2 - (2r-1)(8r^2+4r-22)t + (2r-1)^2 p(r),
\label{eq:abc57} \\
b &=& (t - (r^2-2r)) c + d,
\nonumber\\
a &=& (t - (r^2-2r))^2 c + 2 (t - (r^2-2r)) d
   + (r^2-1)^4  ((4r+4)t + p(r));
\nonumber
\eea
Then the surface is
\be
Y^2 = X^3 + a  X^2 + 8 (r-1)^4 (r+1)^5 b t^2 X + 16 (r-1)^8 (r+1)^{10} c t^4,
\label{eq:S57}
\ee
with a section of height $19/12$ at
\be
X = -\frac{4 (r-1)^4 (r+1)^5 (2r-1) t^2} {(r^2-r+1)^2}
+ \frac{4(r-2)(r+1)^4 t^3}{r^2-r+1} \;.
\label{eq:X57}
\ee
The components of the $A_{11}$ fiber are rational because
the leading coefficient of~$a$ is $9$, a square;
the constant coefficient is $(r^2-r+1)^4 p(r)$,
so $\X(57) / \ang{w_{57}}$ is obtained by extracting
a square root of $p(r)$.  This gives the elliptic curve
with coefficients $[a_1,a_2,a_3,a_4,a_6] = [0,-1,1,-2,2]$,
whose conductor is $57$ as expected
(see e.g.\ Cremona's tables \cite{Cremona}
where this curve appears as \hbox{57-A1(E)}).

This curve has rank~$1$, with generator $P = (2,1)$.
The point at infinity is the CM point of discriminant~$-19$;
this may be seen by substituting $1/s$ for $r$ and
$(t/s^3,X/s^{12},Y/s^{18})$ for $(t,X,Y)$,
then letting $s \ra 0$ to obtain the surface
\be
Y^2 = X^3 + (9t^4-16t^3+4t)X^2 + (72t^5-128t^4)X + (144t^6-256t^5)
\label{eq:X57_19}
\ee
with a $D_6$ fiber at $t=0$ rather than an $A_5$.
Then we still have a section $(X,Y) = (4t^3 - 8t^2, (3t-5)(t^4-t^3))$
of height $19/12$, but there is a \hbox{$2$-torsion} point
$(X,Y)=(-4t,0)$ so $\disc(\Ness) = -\disc(\NS(S))$ drops to
$4 \cdot 12 \cdot (19/12) / 2^2 = 19$.
The remaining rational CM points on $\X(57) / \ang{w_{57}}$
come in six pairs $\pm nP$:
$$
\begin{array}{c||c|c|r|c|c|c}
 n \,&\, 1 \,&\, 2 \,&\,  3 \,&\,  4 \,&\,  5  \,&\,   8 \\ \hline
 r \,&\, 2 \,&\, 1 \,&\, -1 \,&\,  0 \,&\, 5/4 \,&\, 13/9 \\ \hline
-D \,&\, 7 \,&\, 4 \,&\, 16 \,&\, 28 \,&\, 43  \,&\,  163
\end{array}
$$
The last three of these have extra sections $X = -4t$, $X = 0$, and
$$X = -28 \cdot 11^3 \left((t^2/3^6) + (415454 t/3^{18})\right)$$
respectively.  At $r=2$, the $A_{11}$ fiber becomes an $A_{12}$
and our generic \MW\ generator becomes divisible by~$3$;
the new generator $(-972t, 26244t^2)$ has height
$4 - (5/6) - (40/13) = 7/78$,
so $\disc\Ness = 6 \cdot 13 \cdot (7/78) = 7$.
At $r=1$, the $A_5$ and $A_{11}$ fibers together with the section
all merge to form a $D_{18}$ fiber: let $r=1+s$
and change $(t,X)$ to $(st-1, -8s^3 X)$, divide by $(-2s)^9$,
and let $s \ra 0$ to obtain the second Shioda--Hall surface
\be
X^3 + (t^3+8t)X^2 - (32t^2+128)X + 256t
\label{eq:CM4}
\ee
with a $D_{18}$ fiber at $t=\infty$ \cite{Shioda_max,Hall}.
At $t=-1$, the reducible fibers again merge, this time
forming an $A_{17}$ while the \MW\ generator's height drops to
$4 - (4\cdot 14/18) = 8/9$, whence $\disc(\Ness) = 16$.

We find four more rational CM values of~$r$ that do not lift
to rational points on $\X(57) / \ang{w_{57}}$, namely
$r=5$, $1/2$, $17/16$, $-7/4$, for discriminants
$-123$, $-24$, $-267$, and $-627 = -11\cdot 57$ respectively.
The first of these again has an $A_{12}$ fiber, this time
with the section of height $4 - (9/6) - (12/13) = 41/26$; the second
has a rational section at $X=0$; in the remaining two cases
we find the extra section by \hbox{$p$-adic} search:
\be
X = -\frac{11^3 3^2}{2^{21} 91^2} t^2 (7840t^2-2037t+3267) 
\label{eq:CM267}
\ee
for $r=17/16$, and
\be
X =
\frac{3^5 11^4 t^2 q(t)}{2^{12} (81920t^3 + 9216t^2 + 23868t + 39339)^2}
\label{eq:CM627}
\ee
for $r=-7/4$, where $q(t)$ is the quintic
\begin{eqnarray*}
& 419430400t^5 + 2846883840t^4 + 17148174336t^3 \qquad \qquad &
\\
& \qquad \qquad {} + 78784560576t^2 + 175272616341t - 12882888. &
\end{eqnarray*}
Using \cite{Arno} we can show that there are no further
rational CM values.

\section{$N=206$: The last curve $\X(N)/w_N$ of genus zero}

{\bf Summary of results.}
Again we take $\Ness$ of rank~$16$ and an infinite cyclic \MW\ group,
here $R = A_2 \oplus A_4 \oplus A_{10}$ with a \MW\ generator
of height $412/165 = 6 - (1\cdot 2/3) - (2\cdot 3)/5 - (2\cdot 9)/11$.
With the reducible fibers placed at $1,0,\infty$ as usual, the choice
of~$R$\/ means $\Delta = t^5 (t-1)^3 \Delta_1$ with $\Delta_1$ of degree
$24 - (3+5+11) = 5$ and $\Delta_1(0),\Delta_1(1)\neq 0$;
the generator must then have $X(t) = X_1(t)/(t-t_0)^2$ for some
sextic $X_1$ and some $t_0 \neq 0,1$, with the corresponding section
passing through a non-identity component of the $A_2$ fiber
and components at distance~$2$ from the identity of $A_4$ and $A_{10}$.
We eventually succeed in parametrizing such surfaces,
finding a rational coordinate on the modular curve $\X(206)/\ang{w_{206}}$.
These elliptic models do not readily exhibit the involution
$w_2 = w_{103}$ on this curve, so we recover this involution from
the fact that it must permute the branch points of the double cover
$\X(206)$ of $\X(206)/\ang{w_{206}}$.  We locate these branch points as
simple zeros of the discriminant of $\Delta_1$.
As expected, there are $20$ (this is the class number of
$\Q(\kern-.1ex\sqrt{-206})$), forming a single Galois orbit.
We find a unique involution of the projective line
$\X(206)/\ang{w_{206}}$ that permutes these zeros.
This involution has two fixed points, so we switch to
a rational coordinate~$r$ on $\X(206)/\ang{w_{206}}$
that makes the involution $r \, \lra {}-\!r$.  Then $r_0 := r^2$
is a rational coordinate on $\X(206)/\ang{w_2,w_{103}}$,
and the $20$ branch points are the roots of $P_{10}(r^2)$
where $P_{10}$ is the \hbox{degree-$10$} polynomial
\bea
P_{10}(r_0) &=&
8 r_0^{10} - 13 r_0^9 - 42 r_0^8 - 331 r_0^7 - 220 r_0^6 + 733 r_0^5
\label{eq:P10}
\\
&& \quad {} + 6646 r_0^4 + 19883 r_0^3 + 28840 r_0^2 + 18224 r_0\0 + 4096.
\nonumber
\eea
As a further check on the computation, $P_{10}$ has
dihedral Galois group, discriminant $-2^{138} 103^7$,
and field discriminant $-2^{12} 103^5$, while
$P_{10}(r^2)$ has discriminant $2^{311} 103^{14}$
and field discriminant $2^{27} 103^{10}$.
We find that $r=0, \pm 1, \pm 2,\infty$ give CM points of discriminants
$D = -4,-19,-163,-8$ respectively; evaluating $P_{10}(r^2)$
at any of these points gives $-D$ times a square, showing that
the Shimura curve $\X(206)$ has the equation
$s^2 = -P_{10}(r^2)$ over~$\Q$.
The curves $\X(206)/\ang{w_2}$, $\X(206)/\ang{w_{103}}$
are then the double covers
$s_0^2 = -P_{10}(r_0)$, ${s'_0}^{\!2} = -r_0 P_{10}(r_0)$
of the \hbox{$r_0$-line} $\X(206)/\ang{w_2,w_{103}}$
(in that order, because $w_{103}$ cannot fix a CM point
of discriminant $-4$ or $-8$).

{\bf Acknowledgements}

I thank Benedict H. Gross, Joseph Harris, John Voight, Abhinav Kumar,
and Matthias Sch\"utt for enlightening discussion and correspondence,
and for several references concerning Shimura curves and K3 surfaces.
I thank M.~Sch\"utt, Jeechul Woo, and the referees for carefully
reading an earlier version of the paper and suggesting many
corrections and improvements.
The symbolic and numerical computations reported here
were carried out using the packages {\sc gp}, {\sc maxima}, and Magma.

\end{document}